\documentclass[12pt]{article}
\usepackage{microtype} \DisableLigatures{encoding = *, family = * }
\usepackage{subeqn,multirow}
\usepackage{graphicx}
\usepackage{ifpdf}
\ifpdf  \fi
\usepackage{amsfonts,amssymb}
\usepackage{natbib}
\bibpunct{(}{)}{;}{a}{,}{,}

\newcommand{\commentout}[1]{}

\newcommand{\ba}{\begin{array}}
        \newcommand{\ea}{\end{array}}
\newcommand{\bc}{\begin{center}}
        \newcommand{\ec}{\end{center}}
\newcommand{\bdm}{\begin{displaymath}}
        \newcommand{\edm}{\end{displaymath}}
\newcommand{\bds} {\begin{description}}
        \newcommand{\eds} {\end{description}}%17Apr01
\newcommand{\ben}{\begin{enumerate}}
        \newcommand{\een}{\end{enumerate}}
\newcommand{\beq}{\begin{equation}}
        \newcommand{\eeq}{\end{equation}}
\newcommand{\bfg} {\begin{figure}}
        \newcommand{\efg} {\end{figure}}%Nov 5,99
\newcommand{\bi} {\begin {itemize}}
        \newcommand{\ei} {\end {itemize}}
\newcommand{\bpp}{\begin{pspicture}}
        \newcommand{\epp}{\end{pspicture}}
\newcommand{\bqn}{\begin{eqnarray}}
        \newcommand{\eqn}{\end{eqnarray}}
\newcommand{\bqs}{\begin{eqnarray*}}
        \newcommand{\eqs}{\end{eqnarray*}}
\newcommand{\bsq}{\begin{subequations}}
        \newcommand{\esq}{\end{subequations}}
\newcommand{\bsl} {\begin{slide}[8.8in,6.7in]}
        \newcommand{\esl} {\end{slide}}
\newcommand{\bss} {\begin{slide*}[9.3in,6.7in]}
        \newcommand{\ess} {\end{slide*}}
\newcommand{\btb} {\begin {table}[h]}
        \newcommand{\etb} {\end {table}}%Nov 10,99
\newcommand{\bvb}{\begin{verbatim}}
        \newcommand{\evb}{\end{verbatim}}

\newcommand {\pd}[2] {{\frac {\partial {#1}} {\partial {#2}}}}

\newcommand{\reff}[1] {{{Figure \ref {#1}}}}
\newcommand{\refe}[1] {{{(\ref {#1})}}}%Nov 5
\newcommand{\reft}[1] {{{\textbf{Table} \ref {#1}}}}
\def\pmb#1{\setbox0=\hbox{$#1$}%
   \kern-.025em\copy0\kern-\wd0
   \kern.05em\copy0\kern-\wd0
   \kern-.025emkaise.0433em\box0 }

\def\eop{{\hfill $\blacksquare$}}%17Apr01

\newtheorem{theorem}{Theorem}[section]%17Apr01
\def\dt     {{\Delta t}}

\def\T {{{\bf T}}}

\title{On the equivalence between continuum and car-following models of traffic flow} %20140120

\author{Wen-Long Jin \thanks{Department of Civil and Environmental Engineering, California Institute for Telecommunications and Information Technology, Institute of Transportation Studies, University of California, Irvine, CA 92697-3600. Tel: 949-824-1672. Fax: 949-824-8385. E-mail: wjin@uci.edu. Author for correspondence}}

\begin {document}
\maketitle

\begin{abstract}
Recently different formulations of the first-order Lighthill-Whitham-Richards (LWR) model have been identified in different coordinates and state variables. However, there exists no systematic method to convert higher-order continuum models into car-following models and vice versa. In this study we propose a simple method to enable systematic conversions between higher-order continuum and car-following models in two steps: equivalent transformations of variables between Eulerian and Lagrangian coordinates, and finite difference approximations of first-order derivatives in Lagrangian coordinates. 
With the method, we derive  higher-order continuum models from a number of well-known car-following models. We also derive car-following models from higher-order continuum models. For general second-order models, we demonstrate that the car-following and continuum formulations share the same fundamental diagram, but the string stability condition of a car-following model is different from the linear stability condition of a continuum model. This study reveals relationships between many existing models and also leads to a number of new models.
\end{abstract}

{\bf Keywords}:   Eulerian coordinates; Lagrangian coordinates; Primary models; Secondary models; First-order models; Higher-order models; Continuum models; Car-following models; Fundamental diagram; Stability.

\section{Introduction}

\bfg\bc
\includegraphics[width=4in]{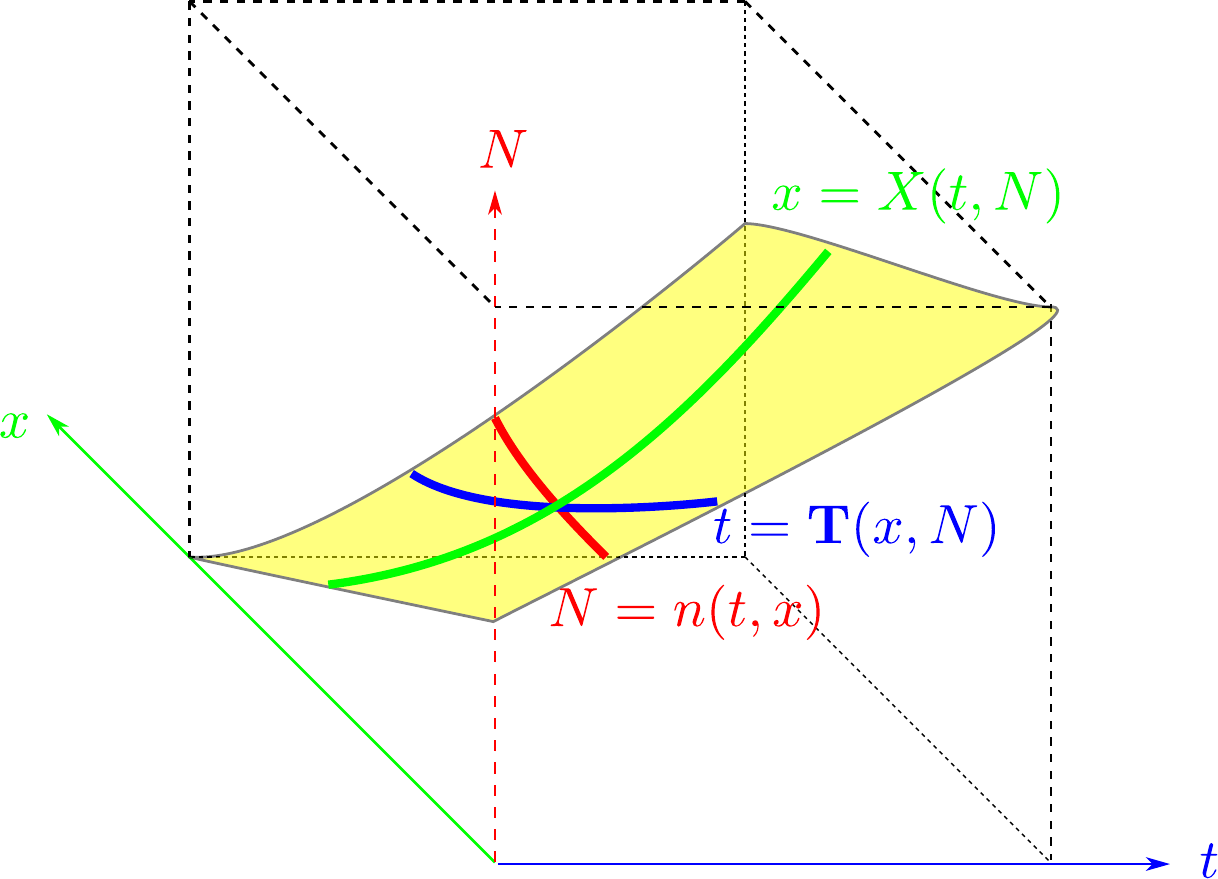}
\caption{The $x-t-N$ surface for a traffic stream} \label{x_t_N_surface}
\ec\efg

In the three-dimensional representation of traffic flow \citep{makigami1971traffic}, as shown in \reff{x_t_N_surface}, the evolution of a traffic stream can be captured by a surface in a three-dimensional space of time $t$, location $x$, whose positive direction is the same as the traffic direction, and vehicle number $N$, which is the cumulative number of vehicles passing location $x$ at $t$ after a reference vehicle \citep{moskowitz1965discussion}. We call these primary variables of traffic flow. In the figure, the red curve, 
\bsq \label{surfaces}
\bqn
N=n(t,x), \label{N-surface}
\eqn
 represents the $x-t$ trajectory of vehicle $N$, which can be collected with GPS devices; the green curve, 
\bqn
 x=X(t,N), \label{x-surface}
\eqn
 represents the cumulative number of vehicles $N$ before $t$ at location $x$, which can be extracted from vehicle counts of a loop detector; and the blue curve,
\bqn 
 t=\T(x,N), \label{t-surface}
\eqn
\esq
  represents the number of vehicles $N$ downstream to $x$ at time $t$, which can be counted from a camera snapshot.
 Among the three primary variables, the two independent variables form a coordinate system and the dependent one is the state variable. In particular, $(t,x)$ form the Eulerian coordinates, and $(t,N)$ the Lagrangian coordinates.

The goal of a traffic flow model is to determine the surface from given initial, boundary, or internal conditions. 
It is usually written as a partial differential equation or a corresponding difference form. 
We call a traffic flow model primary if one of the three primary variables is the unknown state variable, and the order of the model is the maximum order of partial derivatives. 
In contrast, if a partial derivative of a primary variable is used as the unknown state variable, we call the model secondary. 
In this study, we will categorize traffic flow models according to their coordinates, state variables, and orders.

In the Eulerian coordinates, we can define traffic density $k=-n_x$  and flow-rate $q=n_t$, where the subscripts mean partial derivatives. Then the celebrated LWR model \citep{lighthill1955lwr,richards1956lwr} is a partial differential equation in $k$:
\bqn
\pd{k}t+\pd{\phi(k)}x&=&0, \label{E-S-1}
\eqn
which is a hyperbolic conservation law \citep{lax1972shock}. Thus, the LWR model is a first-order secondary model in the Eulerian coordinates, abbreviated as the E-S-1 model.  The model can be derived from the following three rules:
\bi
\item Fundamental law of traffic flow: speed $v=\frac qk$;
\item Conservation of traffic flow: $k_t+q_x=0$;
\item Fundamental diagram of flow-density relation:
\bqn
q&=&\phi(k). \label{flow-density}
\eqn
\ei
In particular, the following triangular fundamental diagram has been shown to match observations in steady traffic well \citep{munjal1971multilane,haberman1977model,newell1993sim}:
\bqs
\phi(k)&=&\min\{v_fk, w(k_j-k)\},
\eqs
where $v_f$ is the free-flow speed, $-w$ the shock wave speed in congested traffic, and $k_j$ the jam density.
The Greenshields fundamental diagram is \citep{greenshields1935capacity}
\bqs
\phi(k)&=&v_f k(1-k/k_j).
\eqs

In \citep{newell1993sim}, the E-P-1 model, i.e., the first-order primary traffic flow model in the Eulerian coordinates, was first discussed:
\bqn
n_t-\phi(-n_x)&=&0, \label{E-P-1}
\eqn
which is a Hamilton-Jacobi equation \citep{evans1998pde}. This model can be numerically solved with the variational principle \citep{daganzo2005variationalKW,daganzo2005variationalKW2}, or the Hopf-Lax formula \citep{claudel2010lax2}.
In \citep{daganzo2006variational}, the L-P-1 formulation of the LWR model, i.e., the first-order primary formulation in the Lagrangian coordinates, was derived based on the inverse function theorem as
\bqn
X_t -\theta(-X_N)&=&0, \label{L-P-1}
\eqn  
which is also a Hamilton-Jacobi equation.
If we denote the spacing by $S=-X_N$, then the speed-spacing relation is
\bqn
\theta(S)=S \phi(1/S).
\eqn
In particular, for the triangular fundamental diagram
\bqs
\theta(S)=\min\{v_f, w(k_jS-1)\};
\eqs
and for the Greenshields fundamental diagram
\bqs
\theta(S)=v_f(1-\frac 1{S k_j}) . 
\eqs
Further in \citep{daganzo2006variational,daganzo2006ca}, various discrete forms of \refe{L-P-1} was obtained based on the variational principle. In particular, Newell's car-following model is shown to be a discrete form of \refe{L-P-1} with the triangular fundamental diagram.
In \citep{leclercq2007lagrangian}, the L-S-1 formulation was derived as a hyperbolic conservation law and numerically solved with the Godunov method \citep{godunov1959}. 
More recently, in \citep{laval2013hamilton}, the equivalent primary and secondary formulations in the $(x,N)$-coordinates, called $T$-models, were derived and solved. In summary, all the primary formulations of the LWR model are Hamilton-Jacobi equations and can be solved with the variational principle, the Hopf-Lax formula, or other techniques \citep{evans1998pde}; and all the secondary formulations are hyperbolic conservation laws and can be solved with the Godunov scheme.

\bfg
\bc
\includegraphics[height=6cm]{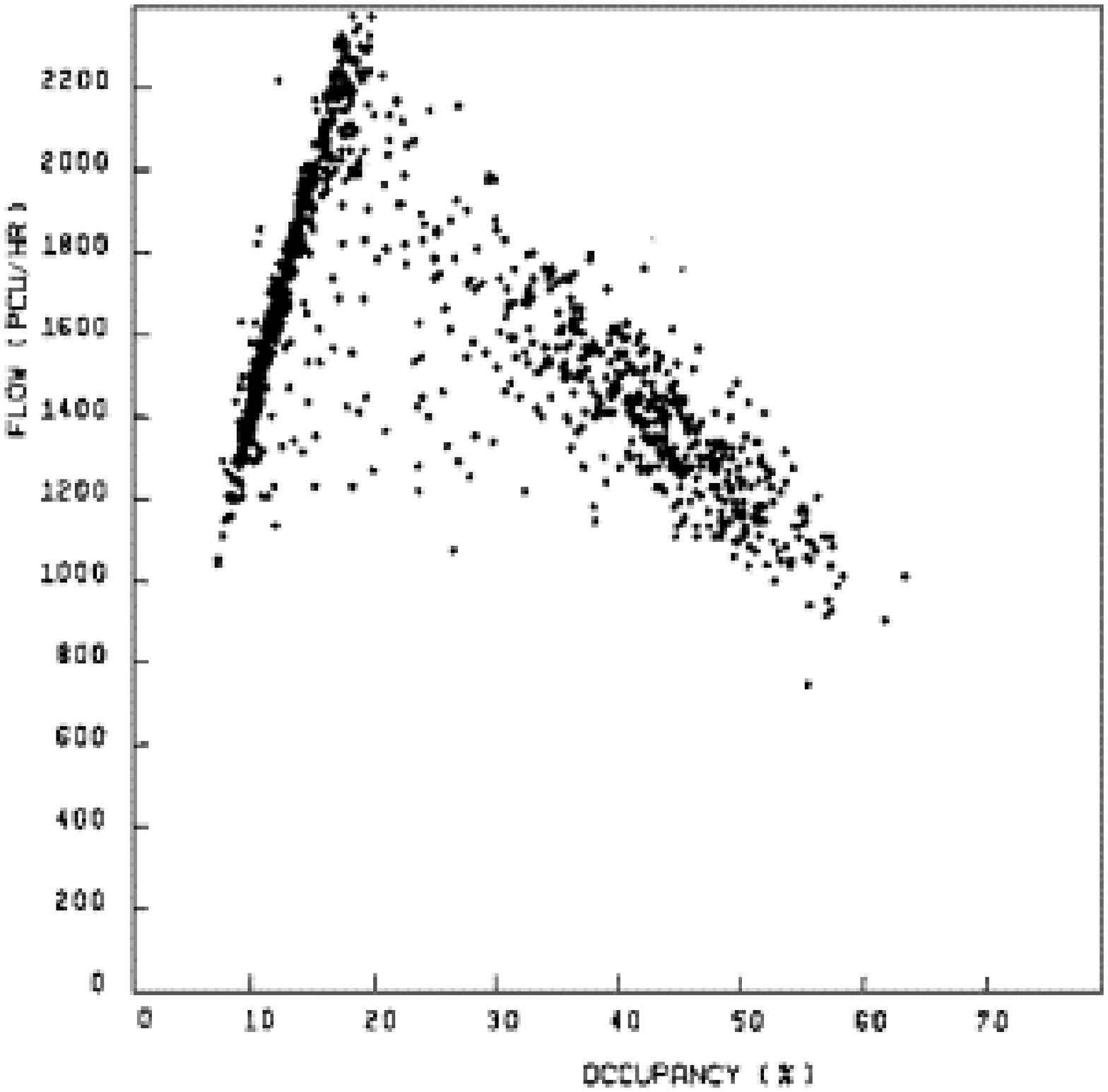}\caption{An observed flow-occupancy relation \citep{hall1986fd}}\label{occupancy_flow}
\ec
\efg

An equilibrium flow-density relationship \refe{flow-density} has been derived from various car-following models in steady states \citep[][Chapter 61]{gazis1961follow,haberman1977model} and verified with empirical observations \citep{delcastillo1995fd_empirical,cassidy1998bivariate}.
However, in reality, traffic can be in non-equilibrium states, and the flow-density or flow-occupancy relations are not unique, as shown in \reff{occupancy_flow}. In such non-equilibrium traffic, instability and hysteresis phenomena have been observed \citep{newell1965instability,treiterer1974hysteresis,sugiyama2008traffic}.  In order to capture such non-equilibrium features of traffic flow, researchers have introduced numerous higher-order continuum models: in \citep[e.g.][]{payne1971PW,whitham1974PW,zhang1998theory,aw2000arz,zhang2002arz}, the conservation equation is coupled with one or more equations in speed or other variables. There have been extensive discussions on theoretical consistency of some higher-order models \citep{daganzo1995requiem,lebacque2007arz,helbing2009controversy,zhang2009comment,helbing2009reply} and on the stability and other properties \citep{kerner1994cluster,li2001stability,jin2003cluster}. Note that, however, these models are still limited, as they are just capable of explaining a subset of many non-equilibrium phenomena \citep{kerner2004physics}. 
Non-equilibrium traffic phenomena have also been successfully modeled with microscopic car-following models since 1950's in the transportation field \citep{pipes1953operational,herman1959traffic} and in the physics field, as traffic systems can be considered many-particle systems \citep{bando1995dmt}. Recently stability and other properties of car-following models have been subjects of some mathematical studies \citep{wilson2001analysis,orosz2006shb}. Reviews of various car-following models are available in  \citep{rothery1992cf,brackstone1999carfollowing,chowdhury2000statistical,helbing2001traffic}.
In the literature, many higher-order continuum models have been derived from car-following models using the method of Taylor series expansion \citep[e.g.][]{payne1971PW,delcastillo1994reaction,zhang1998theory,berg2000continuum,zhang2002arz,jiang2002JWZ}. However, the derivation process is model specific, and different approximation methods can lead to different continuum models for the same car-following model. In addition, such a method cannot be used to derive car-following models from continuum models.

From the studies on different formulations of the LWR model in \citep{daganzo2006ca,daganzo2006variational,leclercq2007lagrangian,laval2013hamilton},
we can see that car-following models can be considered as discrete primary formulations in the Lagrangian coordinates, and continuum models as secondary formulations in the Eulerian coordinates. Based on this observation, in this study we will develop a systematic method for conversions between higher-order continuum models and car-following models. 
However, the framework is different from that for the first-order models: first, the higher-order L-P models are no longer Hamilton-Jacobi equations, and the variational principle or the Hopf-Lax formula cannot be applied to obtain car-following models from the L-P models; second, as higher-order derivatives of the primary variables are involved, the inverse function theorem cannot be used to transform variables between Eulerian and Lagrangian coordinates. 
To resolve the first difficulty, we use simple finite difference methods to derive car-following models from L-P models, since a driver's car-following behavior is only impacted by its leader; to resolve the second difficulty, we systematically derive the transformations of related variables between Eulerian and Lagrangian coordinates.
The conversion method consists of two components: finite difference approximations of L-P derivatives, and transformations between E-S and L-P variables. 
With this method, we will derive continuum formulations of car-following models and vice versa. 
We will be able to reveal relationships between some high-order models and also derive new continuum and car-following models from existing ones.
We will also examine the equivalence between general second-order continuum and car-following models in their steady-state speed-density relations and stability of steady states.

The rest of the paper is organized as follows. In Section 2, we present a method for conversions between continuum and car-following models and apply  it to the first-order models. In Section 3, we derive continuum formulations of some car-following models. In Section 4, we derive car-following formulations of higher-order continuum models. In Section 5, we examine the equivalence between continuum and car-following models. In Section 6, we discuss some future directions.

\section{A systematic conversion method}
In this section we present a systematic method for conversions between car-following and continuum models.
First, we present transformations between E-S and L-P variables, which lead to equivalent E-S and L-P formulations of a traffic flow model. Second, we present finite difference approximations of L-P. derivatives

\subsection{Transformations between E-S and L-P variables}
In  E-S models, variables include $k$, $q$, $v$, and their derivatives with respect to $x$ and $t$; in L-P models, variables include $X$ and its derivatives with respect to $N$ and $t$. In this subsection, we establish the equivalence between E-S and L-P models through transformations between Eulerian and Lagrangian variables.
 
 For a physical variable, we denote it by $F(t,N)$ in the Lagrangian coordinates and by $f(t,x)$ in the Eulerian coordinates. 
  Then from \refe{N-surface} we have
\bqn
f(t,x)&=&F(t,N)=F(t,n(t,x)),
\eqn
which leads to 
\bsq\label{EL-transform1}
\bqn
f_t&=&F_t+F_N n_t,\\
f_x&=&F_N n_x.
\eqn
\esq
In particular, if we let $f(t,x)=x$ and $F(t,N)=X(t,N)$, then from \refe{EL-transform1} we have
\bqs
0&=&X_t+X_N n_t,\\
1&=&X_N n_x,
\eqs
which lead to
\bsq\label{firstderivative}
\bqn
q&=&n_t=-\frac{X_t}{X_N},\\
-k&=&n_x=\frac 1{X_N}.
\eqn
\esq
Substituting \refe{firstderivative} into \refe{EL-transform1}, we obtain the transformation of derivatives in the Eulerian coordinates into those in the Lagrangian coordinates:
\bsq \label{E2L-derivative}
\bqn
f_t&=&F_t- F_N \frac{X_t}{X_N},\\
f_x&=&F_N \frac 1{X_N}.
\eqn
We can repeatedly use \refe{E2L-derivative} to transform higher-order derivatives from E-coordinate to L-coordinate. For example, if we denote $g=f_t$ and $G=F_t- F_N \frac{X_t}{X_N}$, then we have
\bqn
f_{tt}&=&g_t=G_t-G_N \frac{X_t}{X_N}=F_{tt}-2F_{Nt} \frac{X_t}{X_N}+F_{NN}(\frac{X_t}{X_N})^2-F_N \frac{X_{tt}X_N^2-X_t^2 X_{NN}}{X_N^3}.
\eqn
Similarly, we have
\bqn
f_{tx}&=&(F_{tN}-F_{NN}\frac{X_t}{X_N}-F_N\frac{X_{tN}X_N-X_t X_{NN}}{X_N^2}) \frac 1{X_N},\\
f_{xx}&=&\frac{F_{NN}X_N-F_N X_{NN}}{X_N^3}.
\eqn
\esq
In particular, when $f(t,x)=n(t,x)$ and $F(t,N)=N$, we have $q_t=n_{tt}=-\frac{X_{tt}X_N^2-X_t^2 X_{NN}}{X_N^3}$, $q_x=-k_t=n_{tx}=-\frac{X_{tN}X_N-X_t X_{NN}}{X_N^3}$, and $-k_x=n_{xx}=-\frac{X_{NN}}{X_N^3}$.

Similarly we can transform a variable and its derivatives from the Lagrangian coordinates to the Eulerian coordinates as follows:
\bsq
\bqn
F(t,N)&=&f(t,X(t,N)),\\
F_t&=&f_t-f_x \frac{n_t}{n_x}=f_t +v f_x,\\
F_N&=&f_x \frac 1{n_x}= -f_x \frac 1 k,
\eqn
and
\bqn
F_{tt}&=&f_{tt}+ 2 v f_{tx} +v^2 f_{xx}+(v_t+vv_x) f_x,\\
F_{tN}&=&- (f_{tx}+v f_{xx}+v_xf_x) \frac 1k,\\
F_{NN}&=&(k f_{xx} -k_x f_x) \frac 1 {k^3}.
\eqn
\esq
In particular, when $F(t,N)=X$ and $f(t,x)=x$, we have $X_t=v$, $X_N=- \frac 1{k}$, $X_{tt}=v_t+v v_x$, $X_{tN}=- \frac {v_x} k$, and $X_{NN}=- \frac {k_x}{k^3}$.

In \reft{L-E}, where $A=V_t$ and $a=v_t+vv_x$,  we summarize transformations between L-P and E-S variables. 
 In the table, variables are divided into three groups: the first group are for $k,q,v$; the second group for simple derivatives of $k$ and $v$ in the Eulerian coordinates; and the third group for simple derivatives of $X$ in the Lagrangian coordinates. The equivalent variables in E-P and L-S models can be easily obtained through changes of variables but are omitted here. 

\btb\bc
\begin{tabular}{|c|c|c|}\hline
Variables & L-P& E-S\\\hline\hline
Density  &  $-\frac 1{X_N}$  &$k$\\\hline
Speed  & $X_t$  & $v$ \\\hline
Flow-rate  & $-\frac{X_t}{X_N}$  & $q$ \\\hline\hline
Speed rate & $X_{tt}-\frac{X_t}{X_N} X_{tN}$ & $v_t$ \\\hline
Speed gradient & $\frac 1{X_N} X_{tN}$ & $v_x$ \\\hline
Speed curvature & $\frac{X_{tNN}X_N-X_{tN} X_{NN}}{X_N^3}$  & $v_{xx}$ \\\hline
Density rate & $\frac{X_{tN}X_N-X_tX_{NN}}{X_N^3}$  & $k_t$\\\hline
Density gradient  & $\frac{X_{NN}}{X_N^3}$  & $k_x$\\\hline\hline
Acceleration rate  & $X_{tt}$  & $v_t+vv_x$ \\\hline
Jerk  & $X_{ttt}$ & $a_t+va_x$ \\\hline
Speed difference  &$X_{tN}$ &$-\frac{v_x}k$\\\hline
Spacing difference   &$X_{NN}$  & $-\frac{k_x}{k^3}$\\\hline
\end{tabular} \caption{Transformations between E-S and L-P variables} \label{L-E}
\ec\etb

\subsection{Finite difference approximations of L-P derivatives}
In traffic flow, information travels in the increasing direction of $N$; that is, the following vehicles have to adjust their behaviors according to the leading vehicles. Therefore, a simple upwind scheme can be used to numerically solve the L-P formulation as follows: first, we approximate the derivatives with respect to vehicle number by finite differences
\bqs
F_N(t,N)&\approx&\frac {F(t,N)-F(t,N-\Delta N)}{\Delta N}, 
\eqs
In particular, when $\Delta N=1$, a vehicle just follows its immediate leader, and  we have
\bqn
F_N(t,N)&\approx& F(t,N)-F(t,N-1). \label{vehicle-discrete}
\eqn
Similarly we can obtain the vehicle-discrete forms of speed and spacing differences
\bqs
X_N(t,N)&\approx & X(t,N)-X(t,N-1),\\
X_{tN}(t,N)&\approx& X_t(t,N)-X_t(t,N-1),\\
X_{NN}(t,N)&\approx& X(t,N)+X(t,N-2)-2 X(t,N-1).
\eqs
Note that the right-hand sides of the approximations are discrete in vehicle number but still continuous in time. Such vehicle-difference L-P variables are used in time-continuous car-following models.

Furthermore we can use the forward difference to replace the derivatives with respect to time in the vehicle-difference forms:
\bqn
F_t(t,N)&\approx&\frac{F(t+\dt,N)-F(t,N)}{\dt}, \label{time-discrete}
\eqn
where $\dt$ should satisfy the CFL condition \citep{courant1928CFL,leclercq2007lagrangian}: 
$\dt \leq \max |\theta'(S)|$ with $\Delta N=1$, which can be derived from the L-S formulations. In particular, for the triangular fundamental diagram, we have $\dt \leq \frac 1{wk_j}$. Intuitively, information cannot pass for more than one vehicle during a time interval; i.e., a leading vehicle cannot impact the behavior of more than one vehicle, its follower, during a time interval.
With both \refe{vehicle-discrete} and \refe{time-discrete} we can obtain time-discrete car-following models from the L-P formulations.

\subsection{A two-way method}
Based on the transformations between E-S and L-P variables, we can then obtain equivalent E-S and L-P formulations when one of them is given.
From the finite difference approximations of L-P derivatives, we can convert L-P formulations into car-following models and vice versa. Therefore, conversions between continuum and car-following models can be done in two steps, whose order depends on the direction of conversion.
Note that the resulting E-S and L-P models are equivalent, but the car-following models are just approximately equivalent to the L-P models, subject to errors between derivatives and finite differences.

Here we  apply the method to derive equivalent time-continuous and time-discrete car-following models of the LWR model, the first-order continuum model, by discretizing \refe{L-P-1}:
\bqn
X_t(t,N)&=&\theta(X(t,N-1)-X(t,N)),\\
X(t+\dt,N)&=&X(t,N)+ \dt \theta(X(t,N-1)-X(t,N)).
\eqn
In particular, for the triangular fundamental diagram, we have
\bqn
X_t(t,N)&=&\min\{v_f, \frac 1\tau (X(t,N-1)-X(t,N)-S_j)\},\label{pipesmodel}\\
X(t+\dt,N)&=&X(t,N)+ \dt \min\{v_f, \frac 1\tau (X(t,N-1)-X(t,N)-S_j)\}, \label{time-discretepipes}
\eqn
where $\tau=\frac {1}{wk_j}$ is the time gap, and $S_j=1/k_j$ the jam spacing. Here the congested part of the time-continuous car-following model, \refe{pipesmodel}, is essentially the same as the original Pipes car-following model, Equation 2.3 in \citep{pipes1953operational}.
If the minimum time $\dt= \tau$ in the time-discrete Pipes model, \refe{time-discretepipes}, then we have
\bqn
X(t+\tau,N)&=&\min\{X(t,N)+\tau v_f, X(t,N-1)-S_j\}, \label{newellcf}
\eqn
which is Newell's car-following model \citep{newell2002carfollowing,daganzo2006ca}. From the Hopf-Lax formula, we can show that the model is accurate subject to the assumption of uniform spacings between two consecutive vehicles \citep{evans1998pde}.

\section{Continuum formulations of higher-order car-following models}
 In this section, we apply the method in the preceding section to derive continuum formulations of higher-order time-continuous car-following models. In the first step, we approximate finite differences by partial derivatives with respect to vehicle numbers; in the second step, we convert L-P models to E-S continuum models. 

\subsection{Second-order models}
In the literature, a large number of time-continuous car-following models have been developed to calculate
a follower's acceleration rate from its speed, spacing, and speed difference and can be written in the following form:
\bqn
X_{tt}(t,N)&=&\Psi (X_t(t,N),\Delta X(t,N), \Delta X_t(t,N)), \label{2ndcf}
\eqn
where vehicle $N$ follows vehicle $N-1$, the spacing $\Delta X(t,N)=X(t,N-1)-X(t,N)$, and the speed difference $\Delta X_t(t,N)=X_t(t,N-1)-X_t(t,N)$. 

In the first step, from \refe{vehicle-discrete},  we can replace the discrete spacing $\Delta X(t,N)$ by the continuous spacing $-X_N$, and the discrete speed difference $\Delta X_t(t,N)$ by the continuous speed difference $-X_{tN}$. 
Then we obtain the corresponding L-P formulation of the car-following model, \refe{2ndcf}:
\bqn
X_{tt}&=&\Psi(X_t,-X_N, -X_{tN}), \label{2ndL-P}
\eqn
which is a second-order model as the highest order of derivatives is two.
In the second step, from \reft{L-E} we obtain the following E-S formulation:
\bsq\label{2ndE-S}
\bqn
v_t+vv_x&=&\Psi (v,\frac 1 k, \frac{v_x}k), 
\eqn
which, coupled with the conservation equation, 
\bqn
k_t+(kv)_x=0, \label{firstconservation}
\eqn
\esq 
forms a second-order continuum model \citep{payne1971PW,whitham1974PW}.

The following are some examples:
\ben
\item The linear General Motors model without delay, which was derived based on the stimulus-response principle \citep{pipes1953operational,gazis1959car}, can be written as 
\bqn
X_{tt}(t,N)&=&\frac1 T \Delta X_t(t,N), \label{lineargm-cf}
\eqn
where $T>0$ is the reaction time.
Thus $\Psi (v,\frac 1 k, \frac{v_x}k)=\frac 1T \frac{v_x}k$, and the second-order continuum model is (we omit the conservation equation \refe{firstconservation} hereafter.)
\bqn
v_t+\left(v-\frac 1T \frac 1k \right ) v_x&=&0. \label{lineargmcon}
\eqn

\item The nonlinear General Motors model without delay \citep{gazis1961follow} can be written as:
\bqs
X_{tt}(t,N)&=&a \frac{X_t^m \Delta X_t(t,N)}{[\Delta X(t,N)]^l},
\eqs
where $m$ and $l$ are natural numbers, and $a>0$.
Thus $\Psi (v,\frac 1 k, \frac{v_x}k)= a v^m k^{l-1} v_x$, and  the second-order continuum model is
\bqn
v_t+\left(v-a v^m k^{l-1} \right) v_x&=&0. 
\eqn

\item The Optimal Velocity Model (OVM) without delay \citep{bando1995dmt} can be written as:
\bqn
X_{tt}(t,N)&=&\frac1 T (\theta(\Delta X(t,N))-X_t(t,N)). \label{ovm-cf}
\eqn
Thus $\Psi (v,\frac 1 k, \frac{v_x}k)=\frac1 T (\theta(\frac 1k)-v)$, and  the second-order continuum model is
\bqn
v_t+vv_x&=&\frac1 T \left(\theta(\frac 1k)-v\right). \label{ovmcon}
\eqn
Note that the model in \citep{newell1961nonlinear}, $X_t(t+T,N)=\theta(\Delta X(t,N))$, is equivalent to the Optimal Velocity Model if we approximate $X_t(t+T,N)$ by its Taylor series expansion: $X_t(t,N)+T X_{tt}(t,N)$.

\item The Generalized Force Model (GFM) \citep{helbing1998generalized} can be written as:
\bqs
X_{tt}(t,N)&=&\frac1 T (\theta(\Delta X(t,N))-X_t(t,N))-\frac{\Delta X_t(t,N)\Theta(\Delta X_t(t,N))}{T'} e^{-[\Delta X(t,N)-(d+\tau X_t(t,N))]/R'},
\eqs
where $T$ is the acceleration time, $T'<T$ the braking time, $d$ the minimal vehicle distance, $\tau$  the safe time headway (time gap), $R'$ the range of the braking interaction, and $\Theta(\cdot)$ the Heaviside function.
Thus $\Psi (v,\frac 1 k, \frac{v_x}k)=\frac1 T (\theta(\frac 1k)-v)-\frac{v_x\Theta(v_x/k)}{T' k} e^{-[1/k-(d+\tau v)]/R'}$, and the second-order continuum model is
\bqn
v_t+\left(v+\frac 1{T'} \frac{\Theta(v_x/k)}{ k} e^{-[1/k-(d+\tau v)]/R'}\right)v_x&=&\frac1 T \left(\theta(\frac 1k)-v\right).
\eqn

\item The Intelligent Driver Model (IDM) \citep{treiber2000congested} can be written as:
\bqs
X_{tt}(t,N)&=&a\left[1-\left(\frac {X_t(t,N)}{v_f}\right)^\delta -\left(\frac{d+\tau  X_t(t,N)+\frac{X_t(t,N) \Delta X_t(t,N)}{2\sqrt{ab}}}{\Delta X(t,N)}\right)^2 \right],
\eqs
where $a$ is the maximum acceleration rate, $\delta$ the acceleration exponent, $b$ the desired deceleration rate, $\tau$ the safe time headway (time gap), and $d$ the  minimal vehicle distance.
Thus $\Psi (v,\frac 1 k, \frac{v_x}k)=a\left[1-(\frac {v}{v_f})^\delta -(\frac k {k_j}+\tau  k v+\frac{v v_x}{2\sqrt{ab}})^2 \right]$, and the second-order continuum model is
\bqn
v_t+vv_x&=&a\left[1-\left(\frac {v}{v_f}\right)^\delta -\left(d k +\tau  k v+\frac{v v_x}{2\sqrt{ab}}\right)^2 \right]. \label{idmcon}
\eqn
 
\item The Full Velocity Difference Model (FVDM) \citep{jiang2001full} is a combination of the linear General Motors model and the Optimal Velocity Model and can be written as:
\bqn
X_{tt}(t,N)&=&\frac 1 T (\theta(\Delta X(t,N))-X_t(t,N))+\lambda \Delta X_t(t,N), \label{fvdm-cf}
\eqn
where $T$ is the reaction time, and $\lambda$ the sensitivity coefficient.
Thus $\Psi (v,\frac 1 k, \frac{v_x}k)=\frac 1T (\theta(\frac 1k)-v)+\lambda \frac{v_x}k$, and the second-order continuum model is
\bqn
v_t+\left(v-\frac {\lambda}k\right) v_x&=&\frac 1 T\left(\theta(\frac 1k)-v\right). \label{fvdmcon}
\eqn
\een

Among the six continuum models, \refe{lineargmcon}-\refe{fvdmcon}, the one for the linear General Motors model, \refe{lineargmcon}, is a special case of the Aw-Rascle-Zhang model \citep{aw2000arz,zhang2002arz}; the one for the Optimal Velocity model, \refe{ovmcon}, appeared in \citep{Phillips1979traffic}; the one for the Full Velocity Difference Model, \refe{fvdmcon}, is a special case of the Aw-Rascle-Greenberg model \citep{aw2000arz,greenberg2001extensions}; but, to the best of our knowledge, the other three models are new.

\subsection{Third-order models}
In some second-order car-following models, there is a delay term in the acceleration rate:
\bqn
X_{tt}(t+T',N)&=&\Psi (X_t(t,N),\Delta X(t,N), \Delta X_t(t,N)).
\eqn
If approximating $X_{tt}(t+T',N)$ by its Taylor series expansion, $X_{tt}(t,N)+T' X_{ttt}(t,N)$, we then obtain a third-order car-following model:
\bqn
X_{ttt}(t,N)&=&\frac 1{T'}[\Psi (X_t(t,N),\Delta X(t,N), \Delta X_t(t,N))-X_{tt}(t,N)].
\eqn
By using the same method in the preceding subsection, we can then derive the following third-order continuum model:
\bsq
\bqn
k_t+(kv)_x&=&0,\\
v_t+vv_x&=&a,\\
a_t+va_x&=&\frac 1{T'}\left[\Psi (v,\frac 1k, \frac {v_x}k)-a\right],
\eqn
where $a$ is the acceleration rate, and $a_t+va_x$ the jerk.
\esq

For examples, we can derive such continuum models for the General Motors model with delay \citep{gazis1959car,gazis1961follow} and the Optimal Velocity Model with delay \citep{bando1998ovm}. 
Further higher-order models can also be derived by using this method.

\section{Car-following formulations of higher-order continuum models}
In this section, we apply the conversion method in Section 2 to derive car-following models from higher-order continuum models. 
\subsection{A review of higher-order continuum models}

We give a brief review of existing second-order continuum models and categorized them into four families. All second-order models consist of the conservation equation, \refe{firstconservation}, and an additional equation in speed. Here we just list the equation in speed and denote the speed-density relation by $\eta(k)=\theta(1/k)$:
\ben
\item The Phillips family: 
\bi
\item \citep{Phillips1979traffic}: $v_t+vv_x=\frac 1T(\eta(k)-v)$; 
\item \citep{ross1988dynamics}:
$v_t+vv_x=\frac{1}{T} (v_f-v)$,
which was shown to be unrealistic in \citep{Newell1989comments}; 
\item \citep{liu1998high}:
$v_t=\frac{1}{T(k)}(\eta(k)-v)$,
in which the relaxation time $T(k)$ is a function of density $k$.
\ei

\item The Payne family: 
\bi
\item \citep{payne1971PW}: $v_t+vv_x-\frac{\eta'(k)}{2T k}k_x=\frac 1T(\eta(k)-v)$;
\item \citep{whitham1974PW}: $v_t+vv_x+\frac{c_0^2}{k}k_x=\frac 1T(\eta(k)-v)$;
\item \citep{Phillips1979traffic}: $v_t+vv_x+\frac{p'(k)}{k}k_x=\frac 1T(\eta(k)-v)$, where $p(k)$ is the traffic pressure function; 
\item \citep{kuhne1987freeway,kunhe1992continuum,kerner1993cluster}: $v_t+vv_x+\frac{c_0^2}{k}k_x=\frac 1T(\eta(k)-v)+\frac{\mu_0}{k}v_{xx}$, which has a viscous term; 
\item \citep{zhang1998theory}: $v_t+vv_x+\frac{(k \eta'(k))^2}{k}k_x=\frac 1T(\eta(k)-v)$.
\ei

\item The Aw-Rascle family: 
\bi
\item \citep{aw2000arz}: 
\bqn
v_t+(v-k p'(k))v_x=0; \label{ar-original}
\eqn 
\item \citep{zhang2002arz}: 
\bqn
v_t+(v+k \eta'(k))v_x=0; \label{arz}
\eqn
\item  \citep{greenberg2001extensions}: $v_t+(v+k \eta'(k))v_x=\frac 1T(\eta(k)-v)$; 
\item \citep{jiang2002JWZ}: 
 \bqn
 v_t+(v-c_0)v_x=\frac 1T(\eta(k)-v); \label{jwz}
 \eqn
 \item  \citep{xue2003continuum}: $v_t+(v+k\frac{t_r}{T(k)} \eta'(k))v_x=\frac{1}{T(k)}(\eta(k)-v)$, where $t_r$ is driver's reaction time.
 \ei

\item The Zhang family: 
\bi
\item \citep{zhang2003viscosity}: $v_t+(v+2\beta c(k))v_x+\frac{c^2(k)}{k}k_x=\frac 1T(\eta(k)-v)+\mu(k) v_{xx}$, where $c(k)=k \eta'(k)$, $\mu(k)=2\beta T c^2(k)$.
\ei
\een

From these models we have the following observations:
\bi
\item Models of the Phillips family are the simplest: those of the Payne-Whitham family have one extra term with the density gradient $k_x$;  those of the Aw-Rascle family have one extra term with the speed gradient $v_x$; and those of the Zhang family have both.
\item Most models have the total derivative in $v$, $v_t+vv_x$, and some just have the partial derivative, $v_t$. 
\item Some models have the relaxation term, $\frac 1T(\eta(k)-v)$, and some do not.
\item Some models have the viscous term $v_{xx}$, but most do not.
\ei

In the literature, third- or even higher-order continuum models have been proposed \citep{Phillips1979traffic,Helbing1995fluid,Helbing1996ns,colombo2002model,colombo2002traffic}. But in this study we focus on the second-order models.

\subsection{Car-following models of second-order continuum models}
Given a continuum model, in the first step, we convert them into L-P formulations by following \reft{L-E}; 
 in the second step, we approximate partial derivatives with respect to vehicle numbers by finite differences. 
In this subsection we only use the car-following models of the Aw-Rascle family as examples. 

Even though the continuum formulations of General Motors and Full Velocity Difference models belong to the Aw-Rascle family, we can obtain more general car-following models of the following general Aw-Rascle model:
\bqn
v_t+(v-k p'(k))v_x=\frac{1}{T(k)} (\eta(k)-v). \label{aw}
\eqn
First, from \reft{L-E} we can have the L-P formulation of \refe{aw}:
\bqn
X_{tt}&=&\frac{1}{T(-1/X_N)} (\theta(-X_N)-X_t)-p'(-1/X_N) \frac{X_{tN}}{X_N^2}.
\eqn
Further replacing $-X_N$ by $\Delta X(t,N)$ and $-X_{tN}$ by $\Delta X_t(t,N)$, we obtain the following car-following model:
\bqn
X_{tt}(t,N)&=&\frac{1}{T(1/\Delta X(t,N))} (\theta(\Delta X(t,N))-X_t(t,N))+p'(1/\Delta X(t,N)) \frac{\Delta X_t(t,N)}{\Delta X(t,N)^2}. \label{awcf}
\eqn
The new car-following model \refe{awcf} combines a generalized Optimal Velocity model (the first term on the right-hand side) and a nonlinear General Motors model (the second term on the right-hand side). This model is also in the form of \refe{2ndcf}.

For example, the car-following formulation of the Aw-Rascle-Zhang model, \refe{arz}, is
\bqn
X_{tt}(t,N)&=&-\eta'(1/\Delta X(t,N)) \frac{\Delta X_t(t,N)}{\Delta X(t,N)^2}, \label{arzcf}
\eqn
which is different from $X_{tt}(t,N)=\frac{\Delta X_t(t,N)}{T(\Delta X(t,N))}$. The latter was used in  \citep{zhang2002arz} to derive \refe{arz} by Taylor series expansions.  
As another example, the car-following, \refe{jwz}, has the following car-following formulation:
\bqn
X_{tt}(t,N)&=&\frac{1}{T} (\theta(\Delta X(t,N))-X_t(t,N))+c_0 \frac{\Delta X_t(t,N)}{\Delta X(t,N)}, \label{jwzcf}
\eqn
which is a combination of the optimal velocity model and a nonlinear General Motors model.
Note that \refe{jwz} was derived with Taylor series expansion from the FVD model, which is a combination of the optimal velocity model and the linear General Motors model.

To the best of our knowledge, both \refe{arzcf} and \refe{jwzcf} and their corresponding continuous L-P formulations are new. Actually from \refe{arzcf}, we can choose different pressure functions, $p(k)$, to obtain other new car-following models. Following the same procedure, we can obtain more new car-following models from continuum models of other families.

\section{Equivalence between continuum and car-following models}
Different from the first-order models, whose analytical properties and numerical errors have been well understood, 
higher-order models are much more challenging to analyze, as there exist no systematic methods for solving higher-order E-S or L-P models. Thus instead of determining the conditions when continuum and car-following models are equivalent, in section we attempt to examine their equivalence with respect to steady-state speed-density relation and stability of steady states.

\subsection{Steady-state speed-density relations}
In this study we adopt the following definitions of steady states: for car-following models, $X_t(t,N)=v$ and $\Delta X(t,N)=s=\frac 1k$; i.e., all vehicles have the same constant speed $v$ and spacing $s$; for continuum models, $k(t,x)=k$, and $v(t,x)=v$; i.e., the density and speed are constant with respect to both time and location.

For a general second-order car-following model in \refe{2ndcf}, we have $X_tt(t,N)=0$, and $\Delta X_t(t,N)=0$ in steady states. Thus by solving the following equation:
\bqn
\Psi(v,\frac 1k,0)&=&0, \label{speed-density-equation}
\eqn
we can obtain the speed-density relation, $v=\eta(k)$.
For the corresponding continuum model in \refe{2ndE-S}, we have $v_t=0$, and $v_x=0$. Thus the speed-density relation also satisfied \refe{speed-density-equation}. That is, the car-following and continuum model share the same speed-density relation and, therefore, the same fundamental diagram.

In the following we examine three types of higher-order models:
\ben
\item Most of the car-following models in Section 3 and most of the continuum models reviewed in Section 4.1 have a simple relaxation term, $\frac 1T(\theta(\Delta X(t,N))-X_t(t,N))$ or $\frac 1T \left(\eta(k)-v\right)$, and their steady-state speed-density relations are simply $v=\eta(k)=\theta(\frac 1k)$. 
\item For the IDM and its continuum counterpart, \refe{idmcon}, no simple relaxation term is given. We have from \refe{speed-density-equation}
\bqs
1-\left(\frac {v}{v_f}\right)^\delta -\left(d k +\tau  k v\right)^2=0,
\eqs
which leads to the following density-speed relation, $k=\eta^{-1}(v)$:
\bqs
k&=&\frac 1{d+\tau v} \left(1-\left(\frac v{v_f}\right)^\delta\right)^{1/2}.
\eqs
Thus we have the following observations: (i) $v=v_f$ if and only if $k=0$; (ii) $k$ decreases in $v$, and $0\leq k \leq \frac 1d$; (iii) $0<v\leq v_f$ when $k>0$. In particular, when $\delta \to \infty$, we have
\bqs
v&=&\min\left\{v_f, \frac 1\tau \left(\frac 1k-d\right)\right\},
\eqs
which corresponds to the triangular fundamental diagram. The speed-density relation is consistent with that derived from the car-following model in \citep{treiber2000congested}.

\item For the original Aw-Rascle model, \refe{ar-original}, the corresponding car-following model is given by 
\bqs
X_{tt}(t,N)&=&p'(1/\Delta X(t,N)) \frac{\Delta X_t(t,N)}{\Delta X(t,N)^2}.
\eqs
for which the steady states occur if and only if $X_t(t,N-1)=X_t(t,N)$, regardless of the spacing. Therefore, the speed-density relation in steady states is not unique. This problem also appears in the General Motors car-following model \citep{nelson1995stream}. It was argued in \citep{gazis1959car} that the steady-state solutions fall on a valid, unique fundamental diagram if the initial conditions satisfy the same fundamental diagram. However, in a road network where vehicles constantly join and leave a traffic stream at merging and diverging junctions as well as on a multi-lane road, the steady-state speed-density relation cannot be maintained in such models. A simple correction to these models is to add a relaxation term in both continuum and car-following formulations. For example, the FVD model, \refe{fvdm-cf}, can be considered a correction to the linear General Motors model, \refe{lineargm-cf}.

\een

\subsection{Stability of steady states}
In this subsection we analyze the stability of the general second-order car-following model, \refe{2ndcf}, and the corresponding continuum model, \refe{2ndE-S}.

For \refe{2ndcf}, we follow \citep{chandler1958tds} to analyze the string stability of a platoon of vehicles and have the following result.
\begin{theorem}
The car-following model, \refe{2ndcf}, is string stable if and only if 
\bqn
\Psi^2_v(v_0,s_0,0)>2\Psi_s(v_0,s_0,0). \label{stringstability}
\eqn
Note that the string stability condition, \refe{stringstability}, is irrelevant to the partial derivative with respect to the speed difference.
\end{theorem}
{\em Proof}.
 Assuming that the leading vehicle's speed $X_t(t,N-1)=v_0+\epsilon_{N-1} e^{i\omega t}$ and the following vehicle's speed $X_t(t,N)=v_0+ \epsilon_N e^{i\omega t}$ are small monochromatic oscillations around a steady state with a speed $v_0$ and spacing $s_0$, we then have $X_{tt}(t,N)=i\omega \epsilon_N e^{i\omega t}$, $X(t,N-1)=v_0t+\frac{\epsilon_{N-1}}{i\omega} e^{i\omega t} +X(0,N-1)$, and $X(t,N)=v_0t+\frac{\epsilon_{N}}{i\omega} e^{i\omega t} +X(0,N)$. Thus, $\Delta X(t,N)=s_0+\frac{\epsilon_{N-1}-\epsilon_N}{i\omega} e^{i\omega t}$, where $s_0=X(0,N-1)-X(0,N)$ is the initial spacing, and $\Delta X_t(t,N)=(\epsilon_{N-1}-\epsilon_N) e^{i\omega t}$. From \refe{2ndcf}, we have
\bqs
i\omega \epsilon_N e^{i\omega t} &=& \Psi(v_0+\epsilon_N e^{i\omega t},s_0+\frac{\epsilon_{N-1}-\epsilon_N}{i\omega} e^{i\omega t}, (\epsilon_{N-1}-\epsilon_N) e^{i\omega t} ).
\eqs
With the Taylor series expansion of the right-hand side, we then have
\bqs
i\omega \epsilon_N e^{i\omega t} &\approx& \Psi(v_0,s_0,0)+\Psi_v (v_0,s_0,0) \epsilon_N e^{i\omega t}+ \Psi_s(v_0,s_0,0) \frac{\epsilon_{N-1}-\epsilon_N}{i\omega} e^{i\omega t}+\Psi_{\Delta v} (v_0,s_0,0)  (\epsilon_{N-1}-\epsilon_N) e^{i\omega t}.
\eqs
Since $v_0$ and $s_0$ satisfy $\Psi(v_0,s_0,0)=0$ in steady states, we  have
\bqs
i\omega \epsilon_N  &\approx& \Psi_v (v_0,s_0,0) \epsilon_N + \Psi_s(v_0,s_0,0) \frac{\epsilon_{N-1}-\epsilon_N}{i\omega} +\Psi_{\Delta v} (v_0,s_0,0)  (\epsilon_{N-1}-\epsilon_N).
\eqs
Thus we have the ratio of $\epsilon_N$ over $\epsilon_{N-1}$ as
\bqs
\frac{\epsilon_N} {\epsilon_{N-1}}&=& \frac{\Psi_s(v_0,s_0,0)+i\omega \Psi_{\Delta v}(v_0,s_0,0)}{-\omega^2-i\omega \Psi_v(v_0,s_0,0)+\Psi_s(v_0,s_0,0) +i\omega \Psi_{\Delta v}(v_0,s_0,0)}.
\eqs
Then the car-following model, \refe{2ndcf}, is string stable if and only if $|\epsilon_N/\epsilon_{N-1}|<1$ for any frequency $\omega$; i.e., the amplitude of oscillations decays along a platoon. Hence the string stability condition can be written as $\Psi^2_v(v_0,s_0,0)+\omega^2>2\Psi_s(v_0,s_0,0)$ for any $\omega$, 
which is equivalent to \refe{stringstability}.
\eop

For example, for the Optimal Velocity model, \refe{ovm-cf}, the string stability condition, \refe{stringstability}, can be simplified as 
\bqs
2 T \theta_s (s_0)<1.
\eqs
 In particular, for the triangular fundamental diagram, $\theta(s_0)=\min\{v_f, \frac 1\tau (s_0-S_j)\}$, and we have the following conclusion: (i) it is always stable in uncongested traffic, since $\theta_s(s_0)=0$; (ii) in congested traffic, $\theta_s(s_0)=\frac 1\tau$, and the stability condition is equivalent to 
 \bqs
 T<\frac 12 \tau;
 \eqs
i.e., the relaxation time is less than half of the time-gap.

For \refe{2ndE-S}, we follow \citep[][Chapter 3]{whitham1974PW} to analyze its linear stability. We first linearize it  about a steady state at $k_0$ and $v_0$, which satisfy $\Psi(v_0, \frac 1{k_0},0)=0$, by assuming $k=k_0+\epsilon$ and $v=v_0+\xi$, where  $\epsilon(x,t)$ and $\xi(x,t)$ are small perturbations. Omitting higher-order terms, we have \refe{2ndE-S}
\bqs
\epsilon_t +v_0 \epsilon_x +k_0 \xi_x &=&0,\\
\xi_t +(v_0 -\Psi_{\Delta v}(v_0, s_0, 0) s_0) \xi_x &=&\Psi_v(v_0, s_0, 0) \xi- \Psi_s(v_0, s_0, 0) s_0^2 \epsilon,
\eqs 
where $s_0=\frac 1{k_0}$.
Taking partial derivatives of the second equation with respect to $x$ and substituting $\xi_x$ from the first equation we obtain
\bqs
\left(\pd{}{t}+(v_0-\Psi_{\Delta v}(v_0, s_0, 0) s_0) \pd{}{x} \right) \left(\pd{}{t}+v_0 \pd{}{x}\right) \epsilon&=&\Psi_v(v_0, s_0, 0)(\epsilon_t+v_0 \epsilon_x)+ \Psi_s(v_0, s_0, 0) s_0 \epsilon_x.
\eqs
With the exponential solution $\epsilon=e^{i(m x- \omega t)}$, where $m$ is a real number, then the linearized system is stable if and only if two solutions of $\omega$ in the following equation have negative imaginary parts:
\bqs
[\omega-(v_0-\Psi_{\Delta v}(v_0, s_0, 0) s_0) m](\omega-v_0 m) + i [\Psi_v(v_0, s_0, 0)(-\omega+v_0 m)+\Psi_s(v_0, s_0, 0) s_0 m]&=&0,
\eqs
which can be re-written as $\omega^2+(2b_1+i 2b_2) \omega +(d_1+i d_2)=0$ with
\bqs
2b_1&=&-(2v_0-\Psi_{\Delta v}(v_0, s_0, 0) s_0) m, \quad 2b_2=-\Psi_v(v_0, s_0, 0),\\
d_1&=&(v_0-\Psi_{\Delta v}(v_0, s_0, 0) s_0) v_0 m^2, \quad d_2=(\Psi_v(v_0, s_0, 0) v_0+\Psi_s(v_0, s_0, 0) s_0)m.
\eqs
According to \citep{abeyaratne2014macroscopic}, both roots of the equation have negative imaginary parts if and only if $b_2>0$ and $4b_1b_2 d_2-4d_1 b_2^2>d_2^2$; i.e., $\Psi_v(v_0, s_0, 0)<0$, and
\bqs
\Psi^2_s(v_0, s_0, 0)  + \Psi_v(v_0, s_0, 0) \Psi_s(v_0, s_0, 0) \Psi_{\Delta v}(v_0, s_0, 0)<0.
\eqs

Thus for the optimal velocity model, the linear stability condition becomes $\Psi^2_s(v_0, s_0, 0)<0$, which is impossible. 
Such an inconsistency was observed in \citep{abeyaratne2014macroscopic}, where a modified version of \refe{2ndE-S}. However, the modification may not apply to general models. 
Such an inconsistency could be interpreted in two ways: first, it is possible that string stability is different from linear stability; second, it is possible that the continuum and car-following models obtained with the conversion method may not be equivalent, due to the upwind difference method for discretization of the L-P models, \refe{vehicle-discrete}. The real underlying reason is subject to future exploration and seem to be quite challenging.

\section{Conclusion}
Inspired by recent studies on equivalent formulations of the LWR model \citep{daganzo2006ca,daganzo2006variational,leclercq2007lagrangian,laval2013hamilton}, in this study we 
presented a unified approach to convert higher-order car-following models into continuum models and vice versa. The conversion method consists of two steps: equivalent transformations between the secondary Eulerian (E-S) formulations and the primary Lagrangian (L-P) formulations, and approximations of L-P derivatives with finite differences. 
In a sense, continuum and car-following models are equivalent, subject to the errors between finite differences and derivatives. 
The conversion method is different from the traditional method based on Taylor series expansion, as it works in both directions.
We used the method to derive continuum models from general second- and third-order car-following models and derive car-following models from second-order continuum models.
We demonstrated that higher-order continuum and car-following models have the same fundamental diagrams. We also derived the string stability conditions for general second-order car-following models and the linear stability conditions for general continuum models. But we found that that they are different. The underlying reason of such inconsistency is subject to future research, and the conditions when continuum and car-following models are equivalent are yet to be determined. This study further highlights substantial challenges associated with higher-order models.

Through this study, we revealed underlying relationships among many existing models. For example, the Optimal Velocity car-following model \citep{bando1995dmt} and the Phillips second-order continuum model \citep{Phillips1979traffic} are equivalent with the proposed conversion method.  
In addition, with the conversion method we derived a number of new continuum models from existing car-following models and new car-following models from existing continuum models. 
Furthermore, as for the LWR model in \citep{laval2013hamilton}, we can extend the method to derive six equivalent continuous formulations of a traffic flow model: the primary and secondary formulations in $(t,x)$-, $(t,N)$-, and $(x,N)$-coordinates. Different formulations can be helpful for solving different initial-boundary value problems. However, we have yet to develop simple and systematic tools to analyze and solve such higher-order models, as traditional tools for Hamilton-Jacobi equations and hyperbolic conservation laws cannot be directly applied. 

In the future, it is possible to extend the conversion method for more complicated traffic flow models, for example, with multiclass drivers and stochastic driving characteristics, with lane-changing, merging, and diverging behaviors.  We can develop hybrid models from different formulations. We can use different data sources  to calibrate a model and predict traffic conditions.  Finally, we can develop traffic control and management strategies with different formulations.

\section*{Acknowledgments}


\begin{thebibliography}{}

\bibitem[Abeyaratne, 2014]{abeyaratne2014macroscopic}
Abeyaratne, R. (2014).
\newblock Macroscopic limits of microscopic models.
\newblock {\em International Journal of Mechanical Engineering Education},
  42(3):185--198.

\bibitem[Aw and Rascle, 2000]{aw2000arz}
Aw, A. and Rascle, M. (2000).
\newblock Resurrection of ``second order" models of traffic flow.
\newblock {\em SIAM Journal on Applied Mathematics}, 60(3):916--938.

\bibitem[Bando et~al., 1998]{bando1998ovm}
Bando, M., Hasebe, K., Nakanishi, K., and Nakayama, A. (1998).
\newblock {Analysis of optimal velocity model with explicit delay}.
\newblock {\em Physical Review E}, 58(5):5429--5435.

\bibitem[Bando et~al., 1995]{bando1995dmt}
Bando, M., Hasebe, K., Nakayama, A., Shibata, A., and Sugiyama, Y. (1995).
\newblock {Dynamical model of traffic congestion and numerical simulation}.
\newblock {\em Physical Review E}, 51(2):1035--1042.

\bibitem[Berg et~al., 2000]{berg2000continuum}
Berg, P., Mason, A., and Woods, A. (2000).
\newblock Continuum approach to car-following models.
\newblock {\em Physical Review E}, 61(2):1056.

\bibitem[Brackstone and McDonald, 1999]{brackstone1999carfollowing}
Brackstone, M. and McDonald, M. (1999).
\newblock {Car-following: a historical review}.
\newblock {\em Transportation Research Part F}, 2(4):181--196.

\bibitem[Cassidy, 1998]{cassidy1998bivariate}
Cassidy, M. (1998).
\newblock {Bivariate relations in nearly stationary highway traffic}.
\newblock {\em Transportation Research Part B}, 32(1):49--59.

\bibitem[Chandler et~al., 1958]{chandler1958tds}
Chandler, R., Herman, R., and Montroll, E. (1958).
\newblock {Traffic Dynamics: Studies in Car Following}.
\newblock {\em Operations Research}, 6(2):165--184.

\bibitem[Chowdhury et~al., 2000]{chowdhury2000statistical}
Chowdhury, D., Santen, L., and Schadschneider, A. (2000).
\newblock Statistical physics of vehicular traffic and some related systems.
\newblock {\em Physics Reports}, 329(4-6):199--329.

\bibitem[Claudel and Bayen, 2010]{claudel2010lax2}
Claudel, C.~G. and Bayen, A.~M. (2010).
\newblock {Lax--Hopf based incorporation of internal boundary conditions into
  Hamilton-Jacobi equation. part II: Computational methods}.
\newblock {\em IEEE Transactions on Automatic Control}, 55(5):1158--1174.

\bibitem[Colombo, 2002a]{colombo2002model}
Colombo, R. (2002a).
\newblock {A 2$\times$ 2 hyperbolic traffic flow model}.
\newblock {\em Mathematical and Computer Modelling}, 35(5-6):683--688.

\bibitem[Colombo, 2002b]{colombo2002traffic}
Colombo, R. (2002b).
\newblock {Hyperbolic phase transitions in traffic flow}.
\newblock {\em SIAM Journal on Applied Mathematics}, pages 708--721.

\bibitem[Courant et~al., 1928]{courant1928CFL}
Courant, R., Friedrichs, K., and Lewy, H. (1928).
\newblock {{\"U}ber die partiellen Differenzengleichungen der mathematischen
  Physik}.
\newblock {\em Mathematische Annalen}, 100(1):32--74.

\bibitem[Daganzo, 1995]{daganzo1995requiem}
Daganzo, C.~F. (1995).
\newblock Requiem for second-order fluid approximations of traffic flow.
\newblock {\em Transportation Research Part B}, 29(4):277--286.

\bibitem[Daganzo, 2005a]{daganzo2005variationalKW}
Daganzo, C.~F. (2005a).
\newblock A variational formulation of kinematic waves: basic theory and
  complex boundary conditions.
\newblock {\em Transportation Research Part B}, 39(2):187--196.

\bibitem[Daganzo, 2005b]{daganzo2005variationalKW2}
Daganzo, C.~F. (2005b).
\newblock A variational formulation of kinematic waves: Solution methods.
\newblock {\em Transportation Research Part B}, 39(10):934--950.

\bibitem[Daganzo, 2006a]{daganzo2006ca}
Daganzo, C.~F. (2006a).
\newblock {In traffic flow, cellular automata= kinematic waves}.
\newblock {\em Transportation Research Part B}, 40(5):396--403.

\bibitem[Daganzo, 2006b]{daganzo2006variational}
Daganzo, C.~F. (2006b).
\newblock {On the variational theory of traffic flow: well-posedness, duality
  and applications}.
\newblock {\em Networks and Heterogeneous Media}, 1(4):601--619.

\bibitem[Del~Castillo et~al., 1994]{delcastillo1994reaction}
Del~Castillo, J., Pintado, P., and Benitez, F. (1994).
\newblock The reaction time of drivers and the stability of traffic flow.
\newblock {\em Transportation Research Part B}, 28(1):35--60.

\bibitem[{Del Castillo} and Benitez, 1995]{delcastillo1995fd_empirical}
{Del Castillo}, J.~M. and Benitez, F.~G. (1995).
\newblock {On the functional form of the speed-density relationship - II:
  Empirical investigation}.
\newblock {\em Transportation Research Part B}, 29(5):391--406.

\bibitem[Evans, 1998]{evans1998pde}
Evans, L. (1998).
\newblock {\em {Partial Differential Equations}}.
\newblock American Mathematical Society.

\bibitem[Gazis et~al., 1959]{gazis1959car}
Gazis, D., Herman, R., and Potts, R. (1959).
\newblock {Car-following theory of steady-state traffic flow}.
\newblock {\em Operations Research}, 7(4):499--505.

\bibitem[Gazis et~al., 1961]{gazis1961follow}
Gazis, D.~C., Herman, R., and Rothery, R.~W. (1961).
\newblock Nonlinear follow-the-leader models of traffic flow.
\newblock {\em Operations Research}, 9(4):545--567.

\bibitem[Godunov, 1959]{godunov1959}
Godunov, S.~K. (1959).
\newblock A difference method for numerical calculations of discontinuous
  solutions of the equations of hydrodynamics.
\newblock {\em Matematicheskii Sbornik}, 47(3):271--306.
\newblock In Russian.

\bibitem[Greenberg, 2001]{greenberg2001extensions}
Greenberg, J. (2001).
\newblock {Extensions and amplifications of a traffic model of Aw and Rascle}.
\newblock {\em SIAM Journal on Applied Mathematics}, 62(3):729--745.

\bibitem[Greenshields, 1935]{greenshields1935capacity}
Greenshields, B.~D. (1935).
\newblock A study of traffic capacity.
\newblock {\em Highway Research Board Proceedings}, 14:448--477.

\bibitem[Haberman, 1977]{haberman1977model}
Haberman, R. (1977).
\newblock {\em Mathematical models}.
\newblock Prentice Hall, Englewood Cliffs, NJ.

\bibitem[Hall et~al., 1986]{hall1986fd}
Hall, F.~L., Allen, B.~L., and Gunter, M.~A. (1986).
\newblock Empirical analysis of freeway flow-density relationships.
\newblock {\em Transportation Research A}, 20:197.

\bibitem[Helbing, 1995]{Helbing1995fluid}
Helbing, D. (1995).
\newblock {Improved fluid-dynamic model for vehicular traffic}.
\newblock {\em Physical Review E}, 51(4):3164--3169.

\bibitem[Helbing, 1996]{Helbing1996ns}
Helbing, D. (1996).
\newblock {Gas-kinetic derivation of Navier-Stokes-like traffic equations}.
\newblock {\em Physical Review E}, 53(3):2366--2381.

\bibitem[Helbing, 2001]{helbing2001traffic}
Helbing, D. (2001).
\newblock {Traffic and related self-driven many-particle systems}.
\newblock {\em Reviews of Modern Physics}, 73(4):1067--1141.

\bibitem[Helbing, 2009]{helbing2009reply}
Helbing, D. (2009).
\newblock Reply to comment on {��On} the controversy around daganzo��s
  requiem for and {Aw-Rascle��s} resurrection of second-order traffic flow
  models�� by {H.M.} zhang.
\newblock {\em The European Physical Journal B - Condensed Matter and Complex
  Systems}, 69(4):569--570.

\bibitem[Helbing and Johansson, 2009]{helbing2009controversy}
Helbing, D. and Johansson, A. (2009).
\newblock {On the controversy around Daganzo's requiem for and Aw-Rascle's
  resurrection of second-order traffic flow models}.
\newblock {\em The European Physical Journal B}.

\bibitem[Helbing and Tilch, 1998]{helbing1998generalized}
Helbing, D. and Tilch, B. (1998).
\newblock Generalized force model of traffic dynamics.
\newblock {\em Physical Review E}, 58(1):133.

\bibitem[Herman et~al., 1959]{herman1959traffic}
Herman, R., Montroll, E., Potts, R., and Rothery, R. (1959).
\newblock {Traffic dynamics: analysis of stability in car following}.
\newblock {\em Operations Research}, 7(1):86--106.

\bibitem[Jiang et~al., 2001]{jiang2001full}
Jiang, R., Wu, Q., and Zhu, Z. (2001).
\newblock {Full velocity difference model for a car-following theory}.
\newblock {\em Physical Review E}, 64(1):17101.

\bibitem[Jiang et~al., 2002]{jiang2002JWZ}
Jiang, R., Wu, Q.-S., and Zhu, Z.-J. (2002).
\newblock A new continuum model for traffic flow and numerical tests.
\newblock {\em Transportation Research Part B}, 36:405--419.

\bibitem[Jin and Zhang, 2003]{jin2003cluster}
Jin, W.-L. and Zhang, H.~M. (2003).
\newblock The formation and structure of vehicle clusters in the payne-whitham
  traffic flow model.
\newblock {\em Transportation Research Part B}, 37(3):207--223.

\bibitem[Kerner, 2004]{kerner2004physics}
Kerner, B. (2004).
\newblock {\em {The Physics of Traffic: Empirical Freeway Pattern Features,
  Engineering Applications, and Theory}}.
\newblock Springer.

\bibitem[Kerner and Konh{\"a}user, 1993]{kerner1993cluster}
Kerner, B. and Konh{\"a}user, P. (1993).
\newblock {Cluster effect in initially homogeneous traffic flow}.
\newblock {\em Physical Review E}, 48(4):2335--2338.

\bibitem[Kerner and Konh\"auser, 1994]{kerner1994cluster}
Kerner, B.~S. and Konh\"auser, P. (1994).
\newblock Structure and parameters of clusters in traffic flow.
\newblock {\em Physical Review E}, 50(1):54--83.

\bibitem[Kuhne and Michalopoulos, 1992]{kunhe1992continuum}
Kuhne, R. and Michalopoulos, P. (1992).
\newblock Continuum flow models.
\newblock In {\em Traffic flow theory: A state-of-the-art report}, chapter~5.

\bibitem[Kuhne, 1987]{kuhne1987freeway}
Kuhne, R.~D. (1987).
\newblock Freeway speed distribution and acceleration noise: Calculations from
  a stochastic continuum theory and comparison with measurements.
\newblock {\em The Proceedings of the Tenth International Symposium on
  Transportation and Traffic Theory}, pages 119--137.

\bibitem[Laval and Leclercq, 2013]{laval2013hamilton}
Laval, J.~A. and Leclercq, L. (2013).
\newblock The hamilton--jacobi partial differential equation and the three
  representations of traffic flow.
\newblock {\em Transportation Research Part B}, 52:17--30.

\bibitem[Lax, 1972]{lax1972shock}
Lax, P.~D. (1972).
\newblock {\em Hyperbolic systems of conservation laws and the mathematical
  theory of shock waves}.
\newblock SIAM, Philadelphia, Pennsylvania.

\bibitem[Lebacque et~al., 2007]{lebacque2007arz}
Lebacque, J.~P., Mammar, S., and Haj-Salem, H. (2007).
\newblock {The Aw-Rascle and Zhang��s model: Vacuum problems, existence and
  regularity of the solutions of the Riemann problem}.
\newblock {\em Transportation Research Part B}, 41(7):710--721.

\bibitem[Leclercq et~al., 2007]{leclercq2007lagrangian}
Leclercq, L., Laval, J., and Chevallier, E. (2007).
\newblock The lagrangian coordinates and what it means for first order traffic
  flow models.
\newblock {\em Transportation and Traffic Theory}, pages 735--754.

\bibitem[Li, 2001]{li2001stability}
Li, T. (2001).
\newblock $l^1$ stability of conservation laws for a traffic flow model.
\newblock {\em Electronic Journal of Differential Equations}, 2001(14):1--18.

\bibitem[Lighthill and Whitham, 1955]{lighthill1955lwr}
Lighthill, M.~J. and Whitham, G.~B. (1955).
\newblock {On kinematic waves: II. A theory of traffic flow on long crowded
  roads}.
\newblock {\em Proceedings of the Royal Society of London A},
  229(1178):317--345.

\bibitem[Liu et~al., 1998]{liu1998high}
Liu, G., Lyrintzis, A., and Michalopoulos, P. (1998).
\newblock {Improved high-order model for freeway traffic flow}.
\newblock {\em Transportation Research Record: Journal of the Transportation
  Research Board}, 1644:37--46.

\bibitem[Makigami et~al., 1971]{makigami1971traffic}
Makigami, Y., Newell, G.~F., and Rothery, R. (1971).
\newblock Three-dimensional representation of traffic flow.
\newblock {\em Transportation Science}, 5(3):302--313.

\bibitem[Moskowitz, 1965]{moskowitz1965discussion}
Moskowitz, K. (1965).
\newblock { Discussion of `freeway level of service as in uenced by volume and
  capacity characteristics' by D.R. Drew and C. J. Keese}.
\newblock {\em Highway Research Record}, 99:43--44.

\bibitem[Munjal et~al., 1971]{munjal1971multilane}
Munjal, P.~K., Hsu, Y.~S., and Lawrence, R.~L. (1971).
\newblock Analysis and validation of lane-drop effects of multilane freeways.
\newblock {\em Transportation Research}, 5(4):257--266.

\bibitem[Nelson, 1995]{nelson1995stream}
Nelson, P. (1995).
\newblock On deterministic developments of traffic stream models.
\newblock {\em Transportation Research Part B: Methodological}, 29(4):297--302.

\bibitem[Newell, 1965]{newell1965instability}
Newell, G. (1965).
\newblock {Instability in dense highway traffic, a review}.
\newblock In {\em Proceedings}, page~73. Organisation for Economic Co-operation
  and Development.

\bibitem[Newell, 1989]{Newell1989comments}
Newell, G. (1989).
\newblock {Comments on traffic dynamics}.
\newblock {\em Transportation Research Part B}, 23:386--389.

\bibitem[Newell, 2002]{newell2002carfollowing}
Newell, G. (2002).
\newblock {A simplified car-following theory: a lower order model}.
\newblock {\em Transportation Research Part B}, 36(3):195--205.

\bibitem[Newell, 1961]{newell1961nonlinear}
Newell, G.~F. (1961).
\newblock Nonlinear effects in the dynamics of car following.
\newblock {\em Operations Research}, 9(2):209--229.

\bibitem[Newell, 1993]{newell1993sim}
Newell, G.~F. (1993).
\newblock {A simplified theory of kinematic waves in highway traffic {I}:
  General theory. {II}: Queuing at freeway bottlenecks. {III}:
  Multi-destination flows}.
\newblock {\em Transportation Research Part B}, 27(4):281--313.

\bibitem[Orosz and St{\'e}p{\'a}n, 2006]{orosz2006shb}
Orosz, G. and St{\'e}p{\'a}n, G. (2006).
\newblock {Subcritical Hopf bifurcations in a car-following model with
  reaction-time delay}.
\newblock {\em Proceedings of the Royal Society A: Mathematical, Physical and
  Engineering Sciences}, 462(2073):2643--2670.

\bibitem[Payne, 1971]{payne1971PW}
Payne, H.~J. (1971).
\newblock Models of freeway traffic and control.
\newblock {\em Simulation Councils Proceedings Series: Mathematical Models of
  Public Systems}, 1(1):51--61.

\bibitem[Phillips, 1979]{Phillips1979traffic}
Phillips, W. (1979).
\newblock {A kinetic model for traffic flow with continuum implications}.
\newblock {\em Transportation Planning and Technology}, 5(3):131--138.

\bibitem[Pipes, 1953]{pipes1953operational}
Pipes, L. (1953).
\newblock {An operational analysis of traffic dynamics}.
\newblock {\em Journal of Applied Physics}, 24(3):274--281.

\bibitem[Richards, 1956]{richards1956lwr}
Richards, P.~I. (1956).
\newblock Shock waves on the highway.
\newblock {\em Operations Research}, 4(1):42--51.

\bibitem[Ross, 1988]{ross1988dynamics}
Ross, P. (1988).
\newblock {Traffic dynamics}.
\newblock {\em Transportation Research}, 22(6):421--435.

\bibitem[Rothery, 1992]{rothery1992cf}
Rothery, R. (1992).
\newblock {Car following models}.
\newblock {\em Traffic Flow Theory}.

\bibitem[Sugiyama et~al., 2008]{sugiyama2008traffic}
Sugiyama, Y., Fukui, M., Kikuchi, M., Hasebe, K., Nakayama, A., Nishinari, K.,
  Tadaki, S., and Yukawa, S. (2008).
\newblock {Traffic jams without bottlenecks�experimental evidence for the
  physical mechanism of the formation of a jam}.
\newblock {\em New Journal of Physics}, 10(3):033001.

\bibitem[Treiber et~al., 2000]{treiber2000congested}
Treiber, M., Hennecke, A., and Helbing, D. (2000).
\newblock {Congested traffic states in empirical observations and microscopic
  simulations}.
\newblock {\em Physical Review E}, 62(2):1805--1824.

\bibitem[Treiterer and Myers, 1974]{treiterer1974hysteresis}
Treiterer, J. and Myers, J. (1974).
\newblock {The hysteresis phenomenon in traffic flow}.
\newblock {\em Proceedings of the Sixth International Symposium on
  Transportation and Traffic Theory}, page~13.

\bibitem[Whitham, 1974]{whitham1974PW}
Whitham, G.~B. (1974).
\newblock {\em Linear and nonlinear waves}.
\newblock John Wiley and Sons, New York.

\bibitem[Wilson, 2001]{wilson2001analysis}
Wilson, R. (2001).
\newblock {An analysis of Gipps's car-following model of highway traffic}.
\newblock {\em IMA Journal of Applied Mathematics}, 66(5):509.

\bibitem[Xue and Dai, 2003]{xue2003continuum}
Xue, Y. and Dai, S. (2003).
\newblock {Continuum traffic model with the consideration of two delay time
  scales}.
\newblock {\em Physical Review E}, 68(6):66123.

\bibitem[Zhang, 2003]{zhang2003viscosity}
Zhang, H.~M. (2003).
\newblock {Driver memory, traffic viscosity and a viscous vehicular traffic
  flow model}.
\newblock {\em Transportation Research Part B}, 37(1):27--41.

\bibitem[Zhang, 2009]{zhang2009comment}
Zhang, H.~M. (2009).
\newblock {Comment on "On the controversy around Daganzo's requiem for and
  Aw-Rascle's resurrection of second-order traffic flow models" by D. Helbing
  and AF Johansson}.
\newblock {\em The European Physical Journal B}.

\bibitem[Zhang, 1998]{zhang1998theory}
Zhang, H.~M. (1998).
\newblock A theory of nonequilibrium traffic flow.
\newblock {\em Transportation Research Part B}, 32(7):485--498.

\bibitem[Zhang, 2002]{zhang2002arz}
Zhang, H.~M. (2002).
\newblock A non-equilibrium traffic model devoid of gas-like behavior.
\newblock {\em Transportation Research Part B}, 36(3):275--290.

\end{thebibliography}
\end {document}